# Using low Lift-to-Drag spacecraft to perform upper atmospheric Aero-Gravity Assisted maneuvers

Jhonathan O. Murcia Piñeros[1], Rodolpho Vilhena de Moraes[2], Antônio F. Bertachini de Almeida Prado[3]

[1]Corresponding author: jhonathan.pineros@unifesp.br
Institute of Science and Technology ICT-UNIFESP, São José dos Campos, (SP), Brazil.
ORCID: 0000-0002-7013-6515

[2]Institute of Science and Technology ICT-UNIFESP, São José dos Campos, (SP), Brazil.rodolpho.vilhena@gmail.com
ORCID: 0000-0003-1289-8332

[3]National Institute for Space Research INPE, (Graduate Division - DIPGR), São José dos Campos, (SP), Brazil.antonio.prado@inpe.br
ORCID: 0000-0002-7966-3231

**Abstract**
　　The Gravity Assisted Maneuver has been applied in lots of space missions, to change the spacecraft heliocentric velocity vector and the geometry of the orbit, after the close approach to a celestial body, saving propellant consumption. It is possible to take advantage of additional forces to improve the maneuver, like the forces generated by the spacecraft-atmosphere interaction and/or propulsion systems; reducing the time of flight and the need for multiple passages around secondary bodies. However, these applications require improvements in critical subsystems, which are necessary to accomplish the mission. In this paper, a few combinations of the Gravity-Assist were classified, including maneuvers with thrust and aerodynamic forces; presenting the advantages and limitations of these variations. There are analyzed the effects of implementing low Lift-to-Drag ratios at high altitudes for Aero-gravity Assist maneuvers, with and without propulsion. The maneuvers were simulated for Venus and Mars, due to their relevance in interplanetary missions, the interest in exploration, and the knowledge about their atmospheres. The Aero-gravity Assist maneuver with low Lift-to-Drag ratios at high altitudes shows an increase of more than 10° in the turn angle for Venus and 2.5° for Mars. The maneuvers increase the energy gains by more than 15% when compared to the Gravity-Assist. From the Technology Readiness Levels, it was observed that the current level of development of the space technology makes feasible the application of Aero-gravity Assisted Maneuvers at high altitudes in short term**.**

**1. Introduction**

An interplanetary spacecraft projected to travel to other planets, natural satellites, asteroids, or the Sun (celestial bodies) requires different orbits to accomplish its mission. The change from one orbit to another is defined as a maneuver, and it requires impulses and energy. In astrodynamics, the total variation of velocity needed to perform the maneuvers is associated with propellant mass and known as the orbit cost or Delta-V budget (Wertz et al. 2011: pp. 251-253). Large values of Delta-V result in expensive missions. One of the goals of orbit selection and design is to reduce the propellant mass and save costs. There are few options to increase the impulses reducing the burning of propellant, like the Swing-by in the proximity of a secondary body; the manipulation of aerodynamic forces in atmospheric flight, the use of tethers, and/or the use of Solar Radiation Pressure. The two first options will be analyzed in this paper.

The Swing-by or Gravity-Assisted Maneuver (GAM) changes the specific orbital energy of a spacecraft around the massive body (the Sun in this case); varying the heliocentric velocity vector in magnitude and direction, after a close approach to a celestial body (Armellin et al. 2006). The mass of the planet and the proximity of the passage, increase the velocity variations. Throughout the history of spaceflight, the GAM has been implemented successfully in multiple missions, close to natural satellites (ex: Luna 3, Zond 3), asteroids (ex: Galileo, Deep Impact), comets (Vega 1 and 2), and planets (Siddiqi, 2018). Due to the interest of this research, Table 1 presents the closest approaches (or GAM) to the major and minor planets, with an emphasis on the planets Venus, Earth, and Mars, and with approaches below their atmospheres (see Table 1). More details about missions performing swing-by are reported in Siddiqi (2018).

The GAM in major planets, like Jupiter and Saturn, contributes significantly to the velocity change during the passage, due to their larger masses (Broucke, 1988; Broucke and Prado, 1993). For this reason, it is possible to perform the maneuver far from the planet, opposite to the cases of the minor planets (see Table 1). In minor planets, like Earth, Venus, and Mars, the contribution due to their masses is lower, but it is possible to explore the use of their atmospheres to increase the effects of the GAM, for example, the addition of the Lift to the net force gives the possibility of adding forces in the direction or opposite to the gravitational attraction of the planet, which may increase or decrease the turning angle of the maneuver, altering the heliocentric velocity, giving more flexibility for the maneuver. These are some of the effects of aerodynamic forces in GAM, presented by Lewis and McRonald (1991) and Anderson et al. (1991). In this case, the maneuver is called Aero-gravity Assisted Maneuver (AGAM), because it is complemented by Lift and Drag. The missions Galileo, Cassini, and Rosetta passed close enough to be perturbed by the atmosphere of the planet (Table 1), without being projected for the AGAM. These trajectories could be affected by small perturbations due to Drag in Free Molecular Flow (FMF).

Table 1. List of planetary Gravity-Assisted Maneuvers showing the closest approaches.

| Launch year – Spacecraft and Reference | Planet | Closest approach (km) | Date and hour (UTC) |
|---|---|---|---|
| 1972 - Pioneer 10 (NASA, 2021a) | Jupiter | 130 354 | 04/12/1973 - 02:26 |
| 1973 - Pioneer 11 (NASA, 2021b) | Saturn | 20 900 | 01/09/1979 – 16:31 |
| 1973 - Mariner 10 (NASA, 2021b) | Mercury | 327 | 16/03/1975 – 22:39 |
| 1977 - Voyager 1 (NASA, 2021b) | Saturn | 126 000 | 12/11/1980 – 23:46 |
| 1977 - Voyager 2 (NASA, 2021b) | Neptune | 4 800 | 25/08/1989 – 03:56 |
| 1989 - Galileo (NASA, 2021b) | Earth | 303 | 08/12/1992 |
| 1990 - Ulysses (NASA, 2021b) | Jupiter | 378 400 | 08/02/1992 – 12:02 |
| 1997 - Cassini (NASA, 2021b) | Venus | 284 | 25/04/1998 |
| 2004 - MESSENGER (Johns Hopkins, 2021) | Mercury | 199.5 | 06/10/2008 – 08:40 |
| 2004 - Rosetta (NASA, 2021b) | Mars | 250 | 25/02/2007 |
| 2018 - BepiColombo (ESA, 2021a) | Venus | 550 | 10/08/2021-13:48 |
| 2020 - Solar Orbiter (ESA, 2021b) | Earth | 460 | 27/11/2021-04:30 |

The main requirement of the spacecraft to perform an AGAM is the body-shape. It has to be projected to reach large values of Lift-to-Drag ratios (L/D) at the lowest altitudes (Anderson et al. 1991); also, could be complemented by the control of the angle of incidence of the spacecraft relative to the flow (Angle of Attack - AOA), changing the Lift and Drag magnitudes (L/D ratio), and the direction of Lift to manipulate the velocity (Armellin et al. 2006). The GAM, on the other hand, is projected from initial conditions to define the resulting orbit after the swing-by, as a function of the state vector at the pericenter of the planet (which is a function of the angle of approach, the pericenter altitude, and velocity) (Broucke, 1988). The application of the AGAM in minor planets, improves the GAM, increasing the curvature of the GAM and, in some cases, giving more energy to the orbit of the spacecraft (Mazzaracchio, 2015).

Another way to improve the GAM is with a propulsion system; directing the propulsive force (Thrust) and controlling the duration, and magnitude of the force as a function of the mass flow of the propellant. This technique was proposed by Prado (1996) and studied for multiple applications later (Ferreira et al. 2015). The resulting maneuver of the GAM with the propulsive force is known as Powered Gravity-Assisted Maneuver (PGAM) (Prado, 1996). Initially, the PGAM was designed for non-continuum impulses (single or multiple instantaneous burns), but it is also possible to implement continuum thrust, as was presented by Qi and Ruiter (2020), showing an increase in curvature and larger variations of energy compared to the GAM. Detailed surveys of the GAM, PGAM, AGAM and Powered Aero-gravity Assisted Maneuver (PAGAM) were presented by Ferreira et al. (2021); Qi and Ruiter (2021).

As presented before, the GAM has been applied successfully in several space missions, and there are projected extensions of this maneuver, like the AGAM and the PGAM to increase its effects, due to the manipulation of aerodynamic and/or propulsive forces. The scientific literature explores different techniques and configurations of these maneuvers. In this paper are discussed their characteristics and feasibility in section 2. Section 3 presents the mathematical model that describes the dynamics of the system, the atmosphere, and the aerodynamic performance of the selected vehicle. The results of simulating multiple maneuvers with variations in L/D ratios are presented and discussed in section 4, and to close the paper, the conclusions.

In summary, the goal of the present research is to analyze the feasibility and quantify the effects of the AGAM and PAGAM at high altitudes with low values of L/D. The analysis is focused on answering the following questions: Which are the extensions of the GAM and its characteristics? Are they technologically reachable? Which is the range of altitude to observe the effects of the AGAM and PAGAM at low L/D ratios? Which are the effects of these maneuvers? What would be an optimal L/D ratio to increase energy changes?

## 2. Maneuvers derived from the Gravity-Assist

The main effect of the Gravity-Assisted Maneuver (GAM) is the variation of the spacecraft's orbital energy. An energy gain means that the spacecraft increases the velocity after the passage; increasing the semi-major axis and eccentricity of its heliocentric orbit, which is ideal to explore exterior planets and/or to escape from the solar system (Broucke, 1988), for example in missions Voyager 1 and 2. The losses of energy are used to get the inverse behavior, searching close approaches to the internal planets and/or the Sun, like in the missions BepiColombo or Solar Orbiter. Some missions require large values of Delta-V, and it is necessary multiple passages (GAMs) to reach the energy required or to increase the turn angle (direction of the velocity vector after the swing-by), meeting the goals of the mission (Jesick, 2015). When the impulse of the GAM is not coplanar, the maneuvers could be used to modify the orbital plane, changing the inclination (Lohar et al. 1996).

The GAM has been applied in interplanetary travel since 1973; increasing the knowledge and technological development in this area, which is useful for future space missions. The idea of using the aerodynamic forces to modify the GAM was motivated by the increase of studies in hypersonic vehicles like the Waveriders, which present large values of L/D at low Angles Of Attack (AOA below 15°); this is a promising technology useful for future space missions. The hypersonic Waveriders are lifting bodies that attach the shockwave along its leading edge to get large values of L/D (Anderson et al. 1991). However, to date, the first prototype is not developed, due to technical limitations that will be discussed in section 2.1.

The better aerodynamic performance of the Waveriders is achieved at lower altitudes values (larger density of the atmosphere), due to the continuum flow which increases the aerodynamic coefficients, to reach larger L/D ratios. Table 2 presents a simplified review of the maximum values of L/D ratios and the lowest approaches analyzed in the scientific literature.

Table 2. List of AGAM and PAGAM from the scientific literature.

| Planet | Minimum allowed approach (km) | $(Abs(L/D))_{Max}$ | Reference |
|---|---|---|---|
| Titan | -- | 0.4 | Lu and Saikia (2020). |
| Mars | 28 | 3.0 | Henning et al. (2014). |
| Mars | 30 | 3.0 | Qi and Ruiter (2021). |
| Mars | 30 | 3.56 | Armellin et al. (2007). |
| Mars | 40 safety | 3.7 | Armellin et al. (2006). |
| Mars | 60 | 5.0 | Sims et al. (1995). |
| Venus & Mars | -- | 5.0 | Randolph and McRonald (1992). |
| Mars | 30 | 5.38 | Anderson et al. (1991). |
| Venus | 105 | 6.0 | Mazzaracchio (2015). |
| Mars | 20 | 6.63 | Anderson et al. (1991). |
| Venus | 80 | 7.0 | Lewis and McRonald (1991; 1992). |
| Mars | 44 | 7.0 | Lohar et al. (1996). |
| Mars | 60 | 9.0 | McRonald and Randolph (1992). |
| Venus | 100 | 10.0 | Sims et al. (1995). |
| Venus & Mars | -- | 10.0 | Bonfiglio et al. (2000). |
| Mars & Venus | -- | 10.0 | Johnson and Longuski (2002). |
| Mars | < 40 | 10.0 | Lavagna et al. (2005). |
| Venus | 76 | 11.36 | Anderson et al. (1991). |
| Venus | 30 | 14.46 | |
| Venus | 100 | 15.0 | Sims et al. (2000). |

From Table 2 it is possible to determine that Lift-to-Drag ratios larger than 6.8 are required at altitudes lower than 80 km, at Venus and Mars. In the case of altitudes above 100 km at Venus, the mean required values of L/D increase, and are larger than 10.0, because it is necessary to compensate for the decrease in density. However, these values of L/D at upper altitudes are unachievable, due to the dispersion of the particles of the flow, which will be presented in section 3.3.

Another option to improve the GAM is the use of propulsion systems. An impulse in the direction of motion with a radial component pointing to the planet increases the turning angle, as discovered by Prado (1996), which studied the implementation of a single impulse maneuver during the swing-by, to develop the Powered Gravity-Assisted Maneuver (PGAM). If the propulsion system is implemented during an AGAM, the trajectory is defined as a Powered Aero-gravity Assisted Maneuver (PAGAM). The PAGAM combines the advantages of using Lift to increase the effects of the GAM and the implementation of propulsive forces to reduce the losses due to Drag, or to maintain the spacecraft altitude in the atmospheric flight, searching to increase the curvature angle (Mazzaracchio, 2015;

Qi and Ruiter 2021). There are few references from the scientific literature which explored these maneuvers, as presented by: Lewis and McRonald, (1991); Murcia and Prado (2021); Prado, (1996); Prado and Broucke, (1995).

The maneuvers derived from the GAM are presented in Table 3, describing their characteristics and differences. Also, new classifications are suggested, which is a distinguished contribution of the present paper. This classification is important to identify the dominant forces during the passage, limitations of the maneuver, and requirements of the spacecraft. Table 3 shows the differences between the AGAM and PAGAM at high altitude (proposed in this paper), from the other maneuvers.

Table 3. Gravity-Assisted derived maneuvers, requirements, and characteristics.

| Maneuver | Requirements | Characteristics |
|---|---|---|
| Gravity-Assisted Maneuver (GAM) (Broucke, 1988) | - Spacecraft in interplanetary flight near to a secondary body, outside the atmosphere (exoatmospheric flight). | - Changes in the heliocentric orbit due to velocity changes.<br>- Control of the resulting orbit by the altitude of the pericenter, angle of approach, and velocity vector of approach. |
| Aero-gravity Assisted Maneuver at Low Altitudes (AGAM–LA) (Lewis and McRonald, 1991; 1992) | - Spacecraft with an aerodynamic shape to obtain the largest values of L/D, larger than 7.0.<br>- In continuum viscous flow<br>- Maneuver in the planetary atmosphere at low altitudes, below 80 km.<br>- Active aerodynamic control.<br>- Thermal protection. | - Effects of viscous Drag.<br>- High values of heat transference, loads, and dynamic pressure.<br>- Control the Lift direction and magnitude of Lift and Drag with the AOA and Bank angle. |
| Aero-gravity Assisted Maneuver at High Altitudes or Upper Atmosphere (AGAM-HA) | - Conventional spacecraft with flat surfaces or spaceplanes like the shuttle.<br>- Low values of L/D, between -2.0 to 2.0 for low AOA (Anderson, 1991; Blanchard, 1986; Weiland, 2014)<br>- Altitudes larger than 80 km or Knudsen number larger than 0.001. | - Low heat transference.<br>- Low effects of viscous Drag.<br>- Control the Lift direction and magnitude of Lift and Drag with the AOA and Bank angle. |
| Powered Gravity-Assisted Maneuver with Single Impulse (PGAM– SI) (Prado, 1996) | - Single or instantaneous shot.<br>- Control in the direction of the propulsion system.<br>- Also known as powered swing-by (PSB) (Qi and Ruiter 2020) | - Increase or reduction of the GAM effects due to the Delta-V contribution of the maneuver. |
| Powered Gravity-Assisted Maneuver Multi Impulsive (PGAM– MI) (McConaghy et al. 2003). | - Multiple instantaneous burns from the propulsion system.<br>- Control in the direction of the propulsion system. | - It is possible to use different positions during the approach to change the resulting trajectory. |
| | | |
| Powered Gravity-Assisted Maneuver Low thrust (PGAM–LT)(Petropoulos and Longuski, 2004) | - Low thrust and continuum propulsion systems, like electrical or gaseous propulsion systems.<br>- Also described as Low-Thrust Gravity-Assist (LTGA) (McConaghy et al. 2003). | - Actuation along the trajectory increases the effect as a function of propulsion time. |
| Powered Gravity-Assisted Maneuver Conventional Continuum Thrust (PGAM – CT) (Qi and Ruiter 2020) | - Few minutes of propulsion.<br>- Continuum thrust and modulated in magnitude. | - Large consumption of propellant mass.<br>- Increase the effects of Delta V, but at the same time, increase the cost of the mission. |
| Powered Aero-gravity Assisted Maneuver Low Altitudes and Continuum Thrusts (PAGAM – LA/CT) (Qi and Ruiter 2021) | - Same as AGAM at low altitudes and PGAM conventional continuum thrusts. | - Reduce the effects of Drag.<br>- Maintain a constant level of flight to take samples.<br>- Increase the curvature angle. |
| Powered Aero-gravity Assisted Maneuver at High Altitudes and Single Impulse (PAGAM – HA/SI) | - Same as AGAM at high altitudes and Impulsive PGAM.<br>- Murcia and Prado (2017) explored in Free Molecular Flow (FMF) and hypothetical large values of L/D. | - Reduce the propellant consumption.<br>- Low values of L/D and gains of energy with the application of the impulse. |
| Drag – Gravity Assisted Maneuver(DGAM) (Prado and Broucke, 1995) | - Low values of L/D, lower or equal than 1.0. The Drag is the dominant force.<br>- Could be implemented High-Drag devices (like solar sails, parachutes, balloons).<br>- High altitude to reduce the risk of collision, and large values of heat transference and dynamic loads. | - Constant for fixed values of the ballistic parameter.<br>- Periodic, when the spacecraft is rotating and changing periodically the Drag coefficient and/or area.<br>- Controlled with Drag or surface devices like spoilers, to increase or reduce the Drag.<br>- Generate energy losses, useful to capture maneuvers.<br>- Rotate the plane of the orbit. |
| Powered Drag Gravity Assisted Maneuver (PDGAM) | - Same as DGAM<br>- Propulsion system | - Increase or reduce the effects of the DGAM, for example, to reduce the velocity during the passage, collected samples or data from the atmosphere, and after the passage, it increases the energy to escape.<br>- Another application is to null the Drag during low altitude passage, resulting in GAM. |

## 2.1. Risk, limitations and Technological maturity

Table 3 presents several classifications for the AGAM and PGAM, describing the characteristics of each maneuver. Due to the complexity of the maneuvers, it is important to describe their feasibility. The required systems to perform the maneuvers (or critical systems for the maneuver) are analyzed to determine the technological maturity according to the Technology Readiness Levels (TRL) scale (Wertz et al. 2011: pp. 13-16). This method is selected due to the lack of specific classification for space maneuvers. Mankins (2009) describes in detail the nine TRLs. The TRL 1, the lowest one, indicates an infeasible maneuver with the current technology available because only basic principles are observed and reported. The TRL 9 is the highest, corresponding to a feasible maneuver, where all the required systems were proven successful in previous space mission(s).

As presented in Table 1, the Gravity-Assisted Maneuver was performed in several missions, then, the systems required for the maneuver are in TRL 9, showing his feasibility since 1973. The main system to perform the PGAM (single impulsive, multi-impulsive, and continuous) is the propulsion system, widely applied in space missions. With the GAM in TRL 9 and the propulsion systems to the PGAM qualified in space flight (TRL 8), it is possible to say that all of the PGAMs are feasible.

In Table 3, the AGAM and the PAGAM were arranged by altitude, in Low and High Altitudes (lower atmospheric and upper atmospheric flight). The systems required, and the specifications for these maneuvers in atmospheric flight are the high-velocity actuators of the control system (Ex: aerodynamic surfaces); the heat shields and/or high-temperature, strength materials for the airframe (which supports the inertial and heat loads of atmospheric flight, the large values of dynamic pressure, and shock waves; preserving the shape for the aerodynamic performance) (Spilker et al. 2019). In continuum hypersonic flow at low altitudes, it is required a sharped edge geometry to reach large values of L/D (Ex: hypersonic Waveriders), however, it is extremely difficult to manufacture them, due to the small thickness (Santos, 2009). For the same reason, the storage volume is limited, and not enough to assemble conventional systems, like propellant tanks in the case of the PAGAM. To reduce the mass, it is necessary to select low-density and strength materials, with low thermal conductivity and high heat resistance, which is a complex task with the available technology for a slim shape (Lavagna et al. 2005). To date, the unmanned experimental hypersonic vehicles with scramjets are the only system similar to the Waveriders (Voland et al. 2006). These hypersonic vehicles are in the stage of prototype and demonstration (TRL 6), and the hypersonic Waveriders for AGAM at low altitudes are in the concept and application formulation stage, as was presented in Tables 2 and 3 (TRL 2). The control systems for the hypersonic aircraft, the algorithms, and simulations for several cases were documented and reported by Rodriguez et al. (2008), without validation in a relevant environment. For this reason, the aerodynamic control system for the AGAM and PAGAM at low altitudes is classified in TRL 2. The propulsion systems for the atmospheric flight were used successfully in extraterrestrial atmospheres, for the descent of Martian rovers, classified in TRL 8, because it was tested in the operational environment. Another fact to explore is the risk of missing the spacecraft during the maneuver, for example, due to the largest values of the surface temperature during the atmospheric flight, reported being over 10 000 K (Lewis and McRonald, 1992). Also, the variations of the atmospheric conditions due to the solar activity and planetary weather, the uncertainty of the gravitational model, and the planetary topology increase the probability of collision for low passages at hypersonic speeds, requiring high precision sensors and increasing the computational cost for autonomous decisions. As was presented, due to low development of the hypersonic control systems, materials, manufacture, risk of low altitude flight, and low level of technological maturity of the Waveriders; the AGAM, and PAGAM at lowest altitudes are infeasible with the actual technology.

The DGAM, PDGAM, AGAM, and PAGAM at the highest altitudes do not require high values of L/D ratios. Lu and Saikia (2018; 2020) show the possible use of bodies with L/D < 2.0 for aero-maneuvers, as the space capsules, previously qualified on flight (TRL 8). At altitudes larger than 80 km, it is possible the implementation of solar panels, solar sails, and/or Drag devices to control the passage; technology qualified in missions around the Earth (TRL 8). Previous missions applied GAMs in the atmosphere of the planet, showing that maneuvers at upper atmospheric altitudes are feasible (see Table 1). In the upper atmosphere, the risk of collision is reduced because the lowest values of density decrease the dynamic pressure, and the convective heat transference. Due to the TRL of the critical systems, like the ones used in the Space Shuttle (TRL 9), the AGAM and PAGAM maneuvers at the upper atmosphere or high altitudes, are feasible in short term.

As was discussed, the Aero-gravity Assisted Maneuver and Powered Aero gravity Assisted Maneuver at High Altitudes are technologically feasible. For this reason, there are selected to be analyzed in this paper. It is important to clarify that the classification presented here and the analysis made allows the identification of the critical systems to be improved, and the areas of research to route the development of AGAM and PAGAM at Low Altitudes.

### 3. Modeling the maneuvers

The Gravity-Assist was wildly studied and modeled from the Circular Restricted Three-Body Problem (CRTBP) (Broucke, 1988; Szebehely, 1967). Results from simulations show that swing-bys in front of the planet with approach angles in the range from 0 to 180°, reduce the energy of the spacecraft's orbit, and swing-bys behind the planet, in the range from 180° to 360°, increase the energy (Broucke, 1988).

The heliocentric orbit of the spacecraft is integrated numerically from the mathematical model of the CRTBP, according to the following assumptions (Broucke, 1988; Szebehely, 1967):

- The system is composed of the Sun, as the massive body, and the planet, as the secondary body. The mass of the spacecraft is negligible.
- The trajectory is planar.
- The distribution of mass of the primaries is spherically uniform.

- The spacecraft is simulated as a mass point without torques.
- External perturbations are not considered, except by the effect of the atmosphere of the planet during the close approach.
- The Canonical Units of the system are used during the numerical integration.
- The numerical integrator selected is the RKF-7/8 with an adaptive step (Fehlberg, 1968).
- The two planets selected are Venus and Mars, with atmospheres modeled as uniform, varying as an exponential function of the altitude (isothermal atmosphere).

The sequence to calculate the effects of the maneuver and the resulting trajectory are:

I. The spacecraft's initial position is the pericenter of the orbit around the planet, with a tangential velocity of 0.5 VU without perturbation of the atmosphere (VU or Velocity Unit. 1 VU is equivalent to the tangential velocity of the planet around the Sun). The approach angles for the pericenter are 90° and 270°, to observe the minimum and maximum Variations of Energy (VOE), as was determined by Broucke (1988).

II. The orbit is propagated backward to calculate the initial energy and the initial state vector. The initial conditions are calculated when the distance to the center of the planet is higher or equal to 0.5 DU, or when the time is equal to -π/2, (1 DU is equivalent to the planet's semi-major axis) (Broucke, 1988). The state vector at the end of the Sphere Of Influence of the planet (SOI), is used to calculate the curvature angle. The SOI is the spherical region where the planetary gravitational influence is larger than the influence of the Sun (Wertz et al., 2011: p. 281).

III. From the state vector before the passage by the atmosphere, the trajectories are propagated for values of -2.0 ≤ L/D ≤ 2.0. The negative sign indicates that the direction of Lift is against the radius from the planet, a positive Lift means that it is in the radial direction.

IV. During the passage, the values of the close approaches are calculated, as well as the duration of the atmospheric flight, the angle of approach, and the velocity. The altitude is monitored to determine the trajectories ending in collisions.

V. The energy of the spacecraft is calculated after the atmospheric passage, to verify the escape conditions.

VI. The simulation ends when the trajectory of the spacecraft reaches a distance from the center of the planet that is higher or equal to 0.5 DU, or the time after the passage is larger than π/2. At this point, it is calculated the energy and the Variation of Energy (VOE), as presented by Broucke (1988) and Broucke and Prado (1993).

*3.1. The Circular Restricted Three-Body Problem methodology*

The dynamical equations of the planar CRTBP (Szebehely, 1967), as a function of the potential $\Omega$ and the aerodynamic perturbation $A_P$, are:

$$\ddot{x} = 2\dot{y} + \Omega_x + A_{Pi} \qquad (1)$$
$$\ddot{y} = -2\dot{x} + \Omega_y + A_{Pj} \qquad (2)$$

The subscripts $x, y$ indicate the partial derivative and the $i, j$ the component in the direction of $x$ and $y$, respectively. The origin of the reference system is in the center of mass of the primaries and the horizontal axis is in the line connecting the primaries with positive orientation in the direction of the planet. The vertical axis is orthogonal to the horizontal axis and the planet moves counter-clockwise around the Sun. The potential is a function of the spacecraft´s position relative to the Sun $(r_1 = \sqrt{(x+\mu)^2 + y^2})$; relative to the planet $(r_2 = \sqrt{(x-1+\mu)^2 + y^2})$, and to the gravitational constant relative to the mass of the planet $(\mu)$. The potential function is:

$$\Omega = \frac{1}{2}(x^2 + y^2) + \frac{(1-\mu)}{r_1} + \frac{\mu}{r_2} \qquad (3)$$

Figure 1 presents the geometry of the maneuver, including the aerodynamic vectors of Lift (L) and Drag (D); the velocity vector before the approach ($V_i$) and after the approach ($V_f$), the approach angle ($\psi$) and the turn angle ($\delta$).

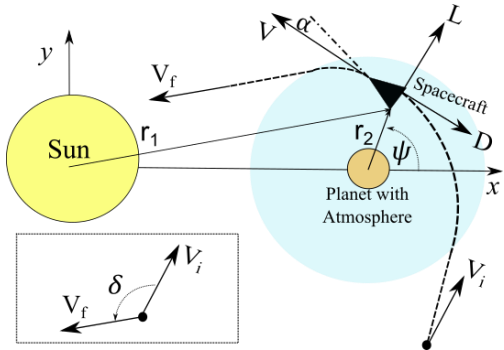

Figure 1. The geometry of the maneuver.

The Variation of Energy per mass unit (VOE) is calculated from the energy before the passage (using the state vector initial conditions), and the energy after the passage, from the state vector at the end of the simulations, according to the methodology presented by Broucke (1988).

For the simulations, the trajectories begin from the pericenter of the orbit around the planet (see list I. to VI.), increasing the altitude to determine the region where the aerodynamic perturbation is affecting the GAM. The energy variations were mapped because they indicate the changes in velocity and geometry of the orbit.

### 3.2. Atmospheric flight

When the spacecraft altitude is below the atmospheric limit of the planet, the acceleration due to aerodynamic forces is significant. The total aerodynamic acceleration is the sum of the acceleration due to Lift force (L/m), in the orthogonal direction relative to the flow, and the acceleration due to Drag force (D/m), in the direction of the incident flow, or against the motion. In this case (m) is the mass of the spacecraft.

$$\vec{A}_P = \frac{\vec{L}}{m} + \frac{\vec{D}}{m} \qquad (4)$$

In Eq. (4), L and D, are functions of the atmospheric density ($\rho$), the aerodynamic coefficients of Lift and Drag ($C_L$ and $C_D$, respectively) which are functions of the Angle Of Attack, (AOA, $\alpha$) and the characteristics of the flow (velocity, molecular distance). The spacecraft Area projected in the direction of the incident flow is A and the magnitude of the flow velocity is $V_\infty$ (Armellin et al. 2006). The aerodynamic relations are:

$$\frac{D}{m} = \frac{1}{2} C_D \left(\frac{A}{m}\right) \rho V_\infty^2 \qquad (5)$$

$$L = D \left(\frac{C_L}{C_D}\right) \qquad (6)$$

The density of the atmosphere is modeled by an exponential function assuming an isothermal atmosphere. In Eq. (7) $H$ is the constant atmospheric scale height, to model the atmosphere of the planet; $h$ is the altitude of the spacecraft from the mean surface of the spherical planet, and $\rho_0$ is the density at the surface (Reagan and Anandarkrishnan, 1993: p.38). Figure 2 presents the density as a function of the altitude for the planets Venus, Earth, and Mars. The values of the scale height and density at the surface were obtained from the NASA Planetary fact sheets (Willians, 2021).

$$\rho = \rho_0 e^{(-h/H)} \qquad (7)$$

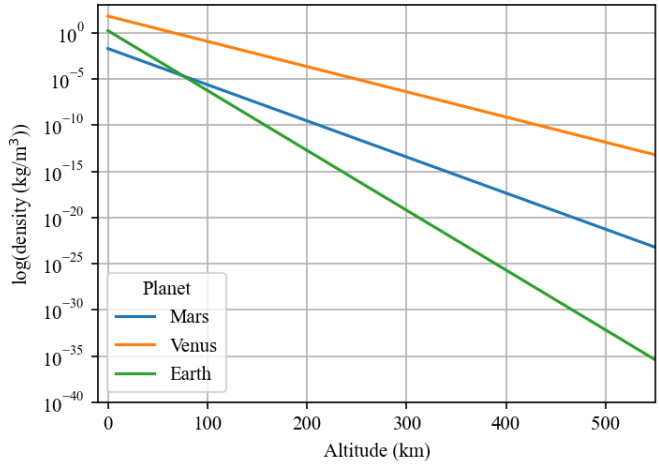

Figure 2. Atmospheric density for the planets Venus, Earth, and Mars.

*3.3. Flow characteristics and aerodynamics*

The simple exponential model of the atmosphere allows to calculate the density as a function of the altitude, showing the largest values of density and aerodynamic accelerations (which are directly proportional) at the lowest altitudes, because the atmosphere has been assumed to be isothermal. Low altitudes reduce the proximity with the planet, increasing the risk of collision. It is also true that the hypersonic velocity in continuum flow increases the loads and heat transfer (which are proportional to density); restricting the mission and the feasibility with the current technology. To reduce the heat transfer, heat shields are usually implemented to protect reentry vehicles, nevertheless, in the case of Waveriders, they represent additional mass and difficulties to the manufacturing process, due to the sharp-edges, as was described by Lewis and McRonald (1992); Santos (2009).

When the spacecraft enters the atmosphere, the average distances of the atmospheric particles are reduced as a function of the increment of the atmospheric density, changing the flow around the vehicle. The general form to classify the flow field is by the Knudsen Number ($K_n$), a non-dimensional parameter which quantifies the dispersion of the particles ($\lambda$) with the characteristic longitude of the spacecraft ($l$); or a relation between the Mach ($M_a$) and the Reynold´s ($R_e$) numbers. The ideal gas heat capacity ratio is $\gamma$, (Reagan and Anandarkrishnan, 1993: p.314; Vivini and Pezzella, 2015: p. 13).

$$K_n = \frac{\lambda}{l} \approx 1.25\sqrt{\gamma}\frac{M_a}{R_e} \qquad (8)$$

The values of $K_n$ classify the fluid in three regimes: Free Molecular Flow (FMF) for values larger than or equal to 10.0 (Prieto et al. 2014; Tewari 2009); transition flow if 10.0 < $K_n$ < 0.01 (Reagan and Anandarkrishnan, 1993: p.310) and hypersonic continuum, when it is lower than or equal to 0.01 (Montenbruck and Gill, 2000: p.84), The $K_n$ is also used to model the aerodynamic coefficients. In general, the FMF is found at the thermosphere and the continuum flow at the stratosphere and troposphere.

The equation to calculate $K_n$ as a function of the exponential density and altitude is available in Reagan and Anandarkrishnan (1993: p.315). The gas molecular weight (Earth = 28.9 gr/mol; Venus and Mars = 43.4 gr/mol for $CO_2$), the effective collision diameter or kinetic diameter (Earth = 3.7 x $10^{-10}$ m; Venus and Mars = 3.3 x $10^{-10}$ m) and a selected reference longitude of the spacecraft of 5 m, are used to calculate the regions for the three flows. Figure 3 presents the Kn as a function of the altitude for the three planets. It is possible to observe that, at Venus, the altitudes of transition regimen are the highest of the three planets, allowing maneuvering the spacecraft farther from the planetary surface.

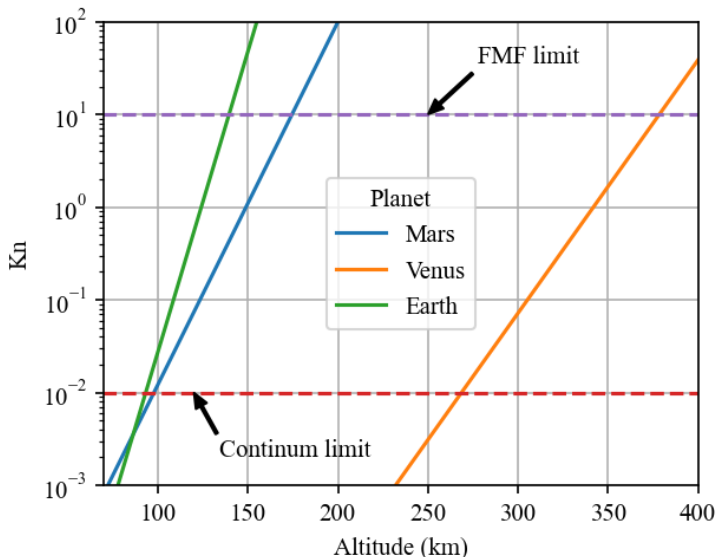

Figure 3. Kn as a function of the altitude for the planets Venus, Earth, and Mars.

In this research, the spacecraft's atmospheric flight is through the three regions of the flow. The approach of the spacecraft begins in the exosphere, at large values of $K_n$, with the flow modeled as FMF. To reduce the losses due to Drag, and because in FMF the L/D ratios are lower than 1.0, it is selected a constant Drag coefficient of 1.0 and without Lift. During the transition flow, the Lift is maintained null and the Drag is modeled using the Bridge formula (Lips and Fritsche, 2005). The last part of the atmospheric flight occurs in the upper atmosphere (high altitude continuum flow), where $10^{-2} \geq K_n \geq 10^{-3}$.

In the continuum flow, the spacecraft coefficients are calculated as a function of the shape of the X-34, based on the Newtonian flow theory presented by Vivini and Pezzella (2015: pp. 89-92). It is observed that the L/D ratio is a function of the Angle of Attack (AOA); reaching values above 2.1, which is a low value compared to the large L/D ratios of the Waveriders, but technologically feasible. Figure 4 presents the values of the aerodynamic coefficients as a function of the AOA, as described by Vivini and Pezzella (2015: pp. 89-92). The L/D for $20° \geq \alpha \geq 0°$, is calculated from:

$$\frac{L}{D} = \frac{C_L}{C_D} = \frac{3.49 \sin^2 \alpha \cos \alpha}{0.046 + 3.49 \sin^3 \alpha} \tag{9}$$

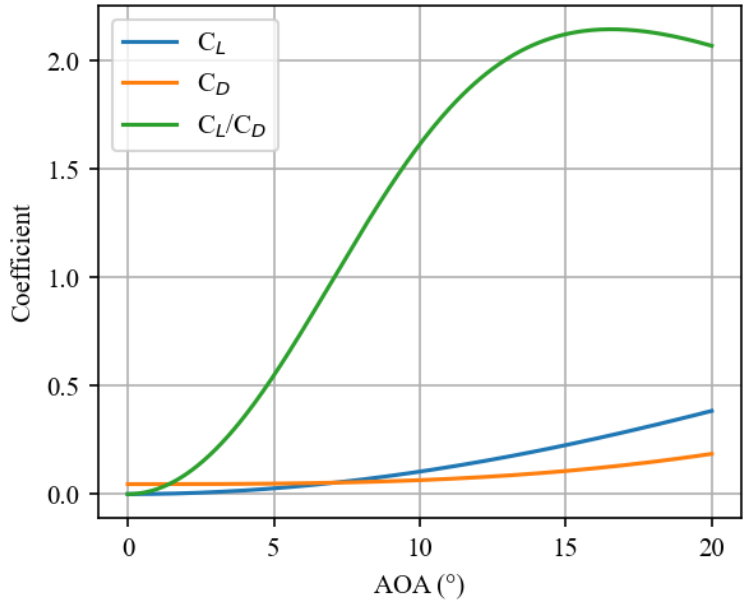

Figure 4. Aerodynamic coefficients.

The first part of this section presents the parameters of the mathematical model for the simulations. An autonomous spacecraft is selected to perform the AGAM and PAGAM. The spacecraft Area-to-Mass ratio (A/m) is 1/50 m²/kg, as presented by McRonald and Randolph (1992). For this paper, the reference length is 5 m. The planets selected to analyze the effects of the AGAM and PAGAM are Venus and Mars, due to the possible applications in future missions. The trajectories are simulated for L/D values between 0.0 to 2.2 (AOA from 0° to 17°). The negative values of L/D ratios indicate an inverted flight (in the direction of the planet) or with a bank angle of 180°, positive L/D is Lift-up or 0° of bank angle. As presented in this Section, during the maneuver in continuum flow at high altitude, the L/D value (also the AOA) and Bank angle are maintained constant to observe the effects under a fixed value of the aerodynamic coefficients along the trajectory. The two approaches angles are 90° and 270°, to verify the application of the methodology developed by Broucke (1988) in the AGAM and PAGAM. To observe the differences between the Aero-gravity Assisted Maneuvers (AGAM) with the Powered Aero-gravity Assisted Maneuvers (PAGAM) at the same altitude, the trajectories were simulated with a continuum thrust in opposite direction and with the same magnitude of the Drag. Then, for the PAGAM, the atmospheric trajectory is dominated by the Lift.

### 4. Effects of the Powered and unpowered Aero-gravity Assisted Maneuvers at High Altitude

In this section are presented and discussed the results from the simulations of the trajectories applying the Aero-gravity Assisted Maneuvers (AGAM) and Powered Aero-gravity Assisted Maneuvers (PAGAM), implementing the models described in section 3.

Figures 5 to 16 show the results of the dependent variable in color scale. The horizontal axis shows the altitude of the approach and the vertical axes the values of the L/D ratios, which are calculated from the Angle of Attack (AOA) using Eq. (9) (see Figure 4). The values of the color scale are presented at the right margin of the figures. For the results was selected the rainbow-scale to observe easily the distribution of the values and regions. The lowest value of the dependent variable is in violet color and the largest in red color.

There are two types of altitudes, the one projected from the Gravity-Assisted Maneuver (GAM), which is the pericenter altitude used to calculate the initial conditions (Section 3 list I. and II.), and the approach altitude, or real altitude, which results from the interaction with the atmosphere. When the Lift is pointing outside the planet (positive L/D or 0° of bank angle), the close approach is larger than the projected altitude, then the maneuver increases the altitude of the approach. For this reason, in the left superior corner of the figures, it is presented a curved white region.

In the first part of this section are presented the results of the maneuvers at Venus, after that the results obtained at Mars. At the end of this section is presented a discussion about the advantages of AGAM and PAGAM over the GAM.

*4.1. Powered and unpowered Aero-gravity Assisted Maneuvers at Venus.*

The first effect observed is the duration of the flight in the continuum flow ($10^{-2} \geq K_n \geq 10^{-3}$) represented as Time of Flight (TOF). This value allows knowing the duration of the interaction between the spacecraft and the particles of the flow, for example, to determine the duration of the influence of the aerodynamic forces or the period where the aerodynamic control is operational. In Figure 5 are presented the TOFs as a function of the L/D ratio and altitude of the approach of the maneuvers, AGAM at the left and PAGAM at the right.

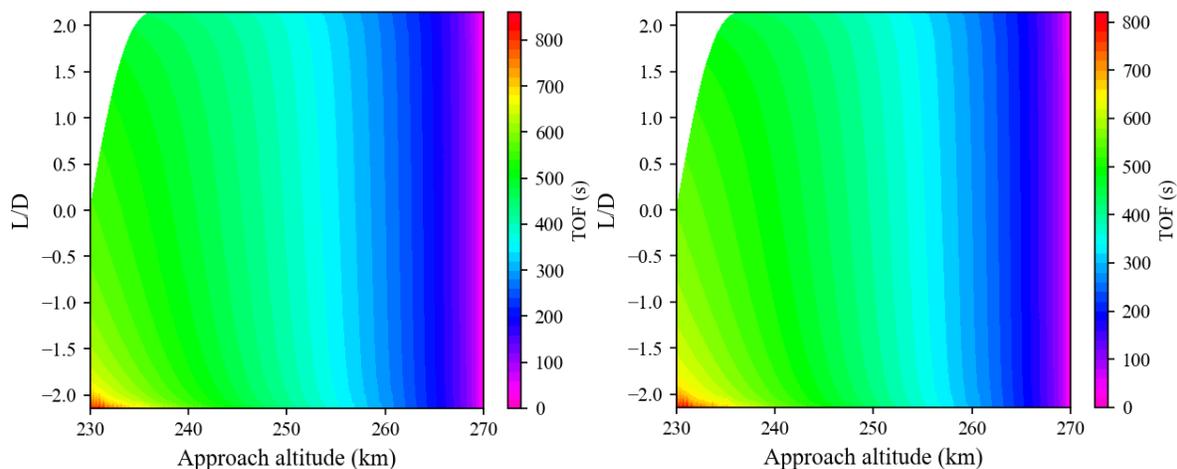

Figure 5. Time of flight for maneuvers approaching at 90° AGAM (left) and PAGAM (right).

Figure 5 shows the largest TOF for the lower value of L/D at the minimum altitude of approach, due to the direction of the Lift being in the same direction of the gravity of the planet, which means that Lift acting in gravity direction increase the atmospheric flight, in double to the positive L/D ratios. The TOF decay with the increase of the L/D ratio, because the Lift in radial direction increases the altitude, reducing the time of the passage.  Above 260 km of altitude, the atmospheric density is lower, and insufficiently to change the TOF as a function of the L/D, then the Lift is insignificant. In the case of the lowest L/D values, the TOF is larger than 800 s, a time that is long enough to take samples from the atmosphere and/or to change the AOA to search optimal trajectories. The only difference between the AGAM and PAGAM is the increase of the TOF for AGAM at the lowest regions of the passage, because the Drag, in the AGAM, reduces the velocity of the passage increasing the TOF, but this is only significant when Lift is maximum and pointing to the planet.

The Variation of Energy (VOE) by the mass unit is an important measurement for the Gravity-Assisted Maneuvers (GAM) because it represents the changes in the heliocentric velocity with the variations of the semi-major axis and eccentricity of the orbit.  In astrodynamics, the increase of orbital energy is equivalent to an increase in the values of the eccentricity, turning the orbit more hyperbolic, for this case. If the energy is reduced, the orbit reduces the eccentricity, semi-major axis, and velocity. Broucke (1988) shows that GAM with an approach angle of 90° gives the maximum loss of energy, more details were presented in Section 3.

The change in the direction of the heliocentric velocity vector is measured by the turn or the curvature angle, which is the angle between the velocity vector before the approach with the vector after the approach. Large values of the curvature angle indicate that the maneuver is useful to direct the spacecraft to other parts of the solar system. Figure 6 shows the VOE, and Figure 7 shows the curvature angles, at 90° of approach at Venus.

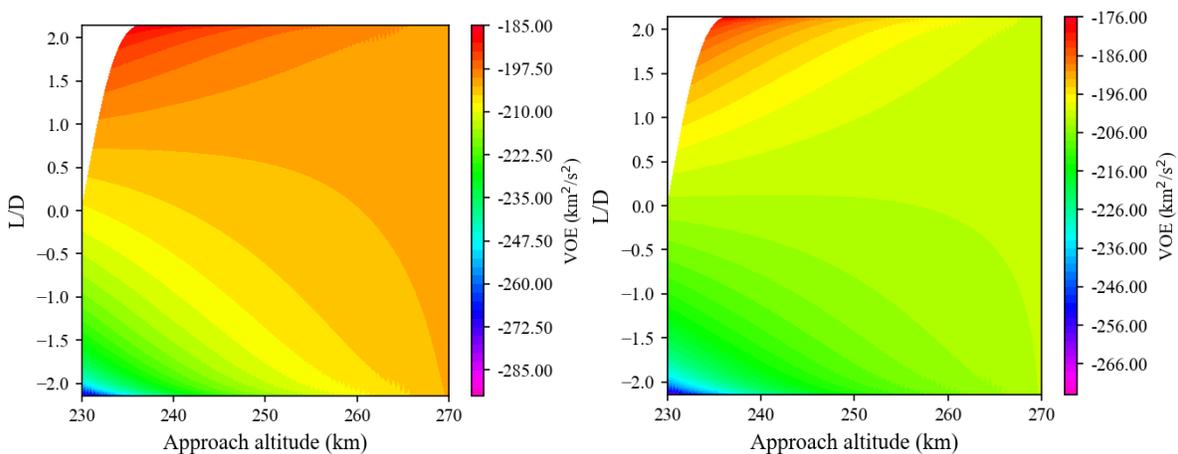

Figure 6. Variation of energy for the angle of approach of 90° for AGAM (left) and PAGAM (right).

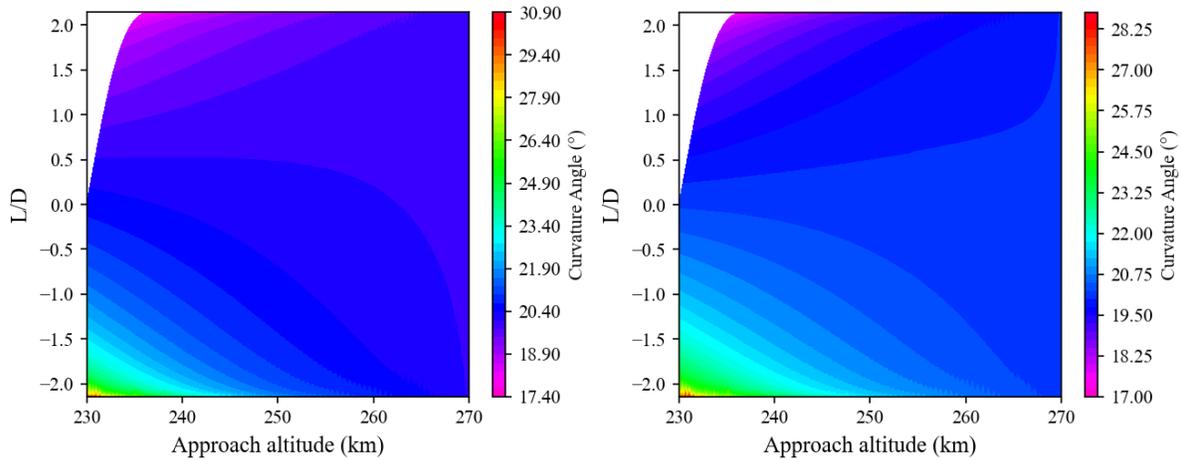

Figure 7. Variation of the angle of curvature when the angle of approach is 90° for AGAM (left) and PAGAM (right).

The AGAM, like the GAM, decreases the VOE after the passage for an angle of approach of 90°, so minimizing the VOE at the lowest altitudes with the lower AOA (or maximum negative L/D ratio). The PAGAM reduces the energy losses because the Drag is null, due to the effect of the Thrust, see Figure 6. In the AGAM, the variations of Energy from 0.8 to 1.1 L/D ratios, are similar to the VOE at larger altitudes, like 270 km, where the influence of Lift decreases. In the case of the PAGAM, the region is reduced, between 0.1 and 0.4 L/D ratios. The losses of Energy increase the curvature angle, and it is maximum in the closest approach with the lowest values of L/D (see Figure 7). As the PAGAM increases the energy, the curvature angle is reduced because the energy gains increase the eccentricity of the orbit. In Figure 7, the blue color is closer to the effect of the GAM. Comparing this region with the red one, the AGAM increases the curvature angle by more than 10° and the PAGAM above 8°, to the GAM, which is an interesting result, because different L/D configurations allow to approaching the maneuvers the GAM, which is a technique useful, for cases of malfunctions on the guidance system and/or atmospheric variations. This region could be a safety zone for the maneuver.

The methodology described for the GAM was applied to select the pericenter, in this case for altitudes from 230 km to 270 km, and approach angle of 90 (Broucke, 1988; Prado, 1996). Due to the influence of the atmosphere, it is possible to change the real approach angle and altitude from the projected one for the maneuver (pure GAM). The differences between the pericenter of the AGAM to the pericenter of the GAM (in altitude and angle of the approach), are presented in Figure 8.

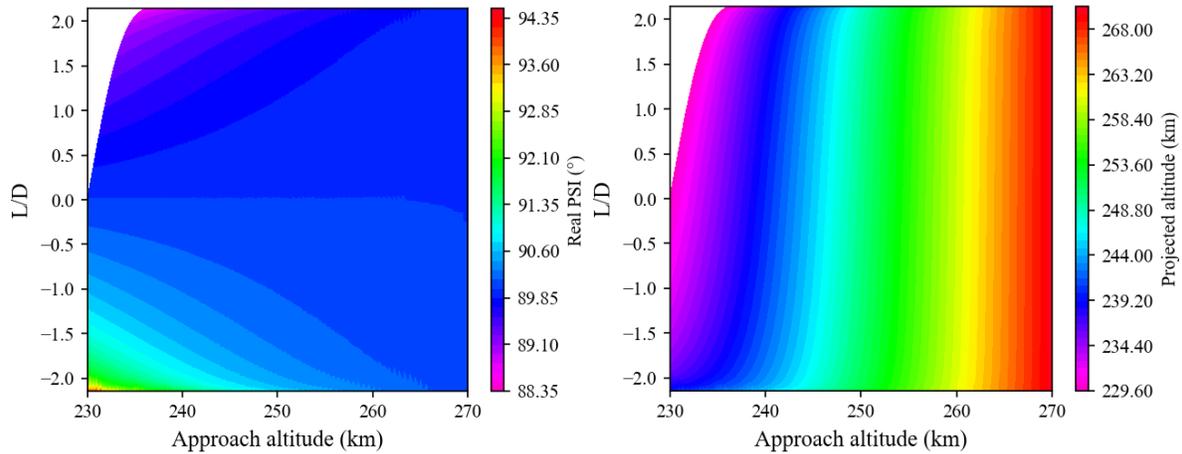

Figure 8. Variation of approach angle (left) and altitude of the pericenter (right) for AGAM at 90°.

Figure 8 shows that the approach angle increases at the closest approach, due to the increase in the TOF and the influence of Lift which increases the atmospheric flight. The lowest values of L/D ratios increase the approach angle by around 4° from the projected value (or GAM), and the largest values reduce the angle of approach by more than 1° (see Figure 8). The variation is lower than 4% from the desired value (90°). The results represent changes in the position of the pericenter, which indicate a small rotation of the orbital plane, due to the lateral maneuvers restricted to the plane of the orbit. In the case of the projected altitude, it is observed a decrease in negative values of L/D, because the direction of the Lift increases the net force in direction of the planet. In the case of a projected altitude of 230 km at the lowest L/D, the approach altitude is below the altitude limit of the flow to be analyzed, and then these values are not considered here. For 240 km of projected altitude, the real approach is below 234 km (see figure 8, right side). In the case of positive values of L/D, the approach increases the altitude, decreasing his influence with the altitude because the Lift is lower due to the density reduction. The same behavior is observed for PAGAM without significant differences.

For an approach angle of 270°, the effects in the Time Of Flight, curvature angle, approach angle, and approach altitude are equivalent to the ones observed at 90° of approach angle. The only significant differences are observed in the Variations Of Energy (VOE). Figure 9 shows the VOE for the approach behind the planet, or 270°, where, according to the GAM, are expected gains of energy.

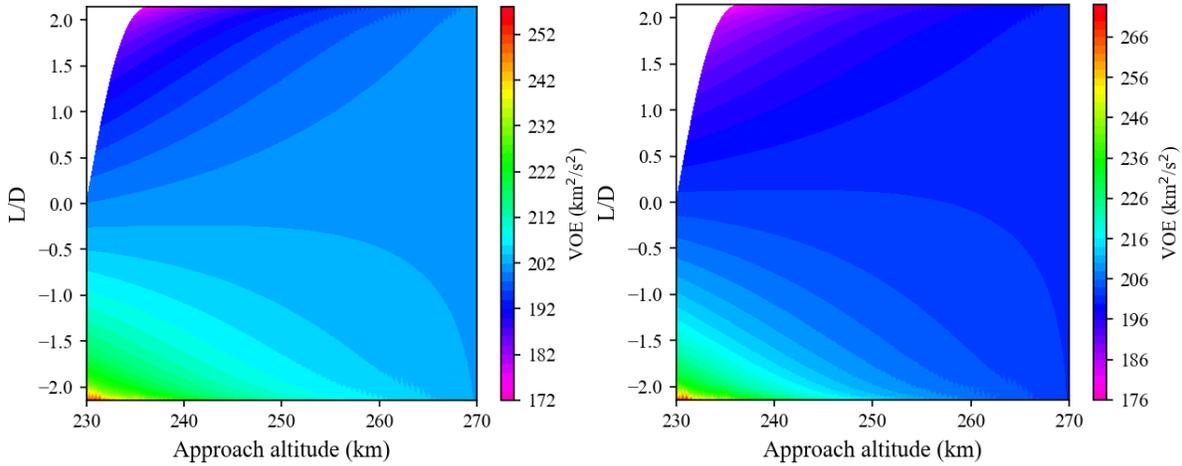

Figure 9. Variations of energy at 270° AGAM (left) and PAGAM (right).

In this approach, the L/D pointing to the planet increases the VOE by more than 50 km²/s², when it is compared to the GAM (largest altitudes or blue region). The PAGAM increases the energy of the AGAM, presenting gains of more than 4 km²/s² for the lowest L/D and more than 14 km²/s² for the largest values (Figure 9), due to the increase of the effect of the net force in the direction of the planet. Then, the PAGAM increases considerably the energy of the AGAM.

*4.2. Powered and unpowered Aero-gravity Assisted Maneuvers at Mars*

As Venus, the planet Mars was selected to observe the effects of the maneuvers. Compared to Venus, Mars has a lower mass, a large distance to the Sun, and different atmospheric properties. The effects of the AGAM and PAGAM for low L/D ratios at the Upper Atmosphere are similar, but lower in magnitude.

The results obtained from the simulations show that the TOF is maximum at 75 km of approach altitude with the lowest L/D ratio (more than 450s), and the effect is reduced above 98 km of altitude. The approach angles change in about 1.5° from the projected value, and the mean curvature angle without effects of Lift (GAM) is around 9°, increasing by more than 2.5° for the minimum L/D and decreasing around 1.2° for maximum. For these reasons, only the effects over the energy are presented in Figure 10, for an angle of approach of 90°, and Figure 11, for an angle of 270°.

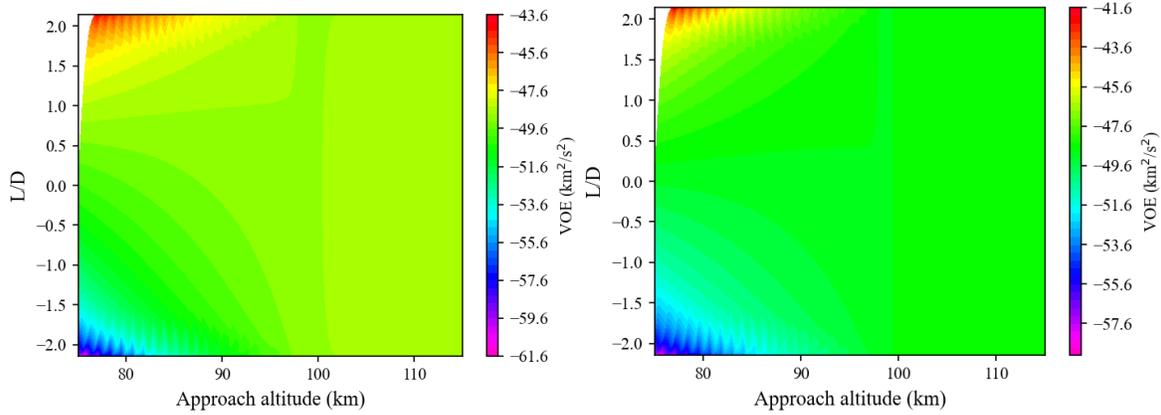

Figure 10. Variations of energy after passage around Mars, angle of approach of 90° for AGAM (left) and PAGAM (right).

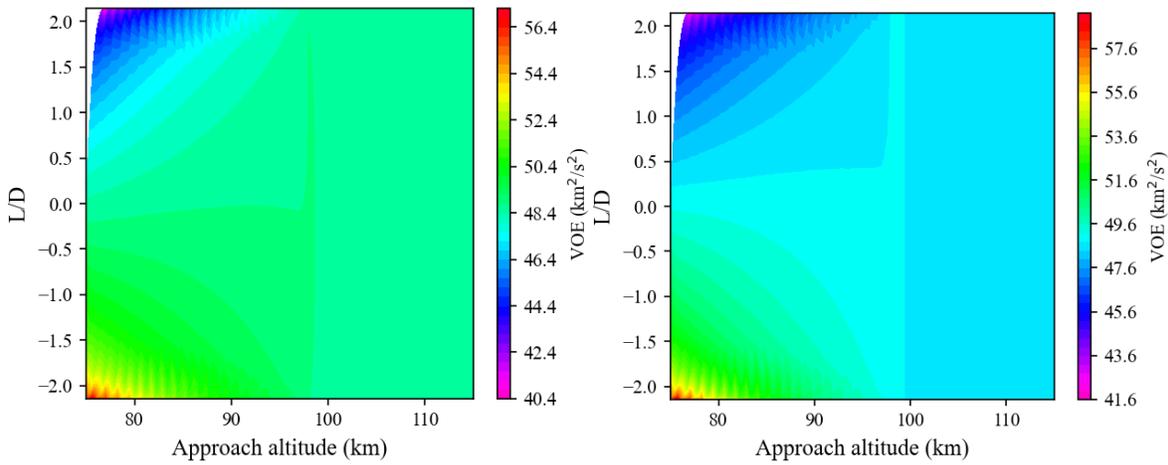

Figure 11. Variations of energy after passage around Mars, angle of approach of 270° for AGAM (left) and PAGAM (right).

At Mars, like at Venus, the PAGAM reduces the losses of energy when is compared to the AGAM, and the maximum and minimum regions are located at the maximum L/D, for 0° and 180° of bank angle, at lowest altitudes, at larger atmospheric density. In this case, the maximum values of L/D don't increase significantly the real altitude from the projected one, and the white region is small (see Figures 10 and 11), also, the effects of the LIft only are observed below 98 km of altitude due to the density of the atmosphere.

### 4.3. Contributions of the Powered and unpowered Aero-gravity Assisted Maneuvers to the Gravity Assisted Maneuver

The changes in VOE and curvature angle due to the PAGAM and the AGAM are compared to the GAM, to quantify the contributions of the maneuvers. The initial conditions of approach are the same for the three maneuvers, like altitude, pericenter velocity, and approach angle. Then, the GAM is analyzed hypothetically without the influence of the atmosphere (or Thrust equal to Drag and without Lift). The results of the contribution in energy are presented in Figures 12, for Venus, and 13, for Mars, in terms of percentage, when compared to the value of the GAM at the same point. Figures 14 and 15 present the contributions of the maneuver in the turn angle.

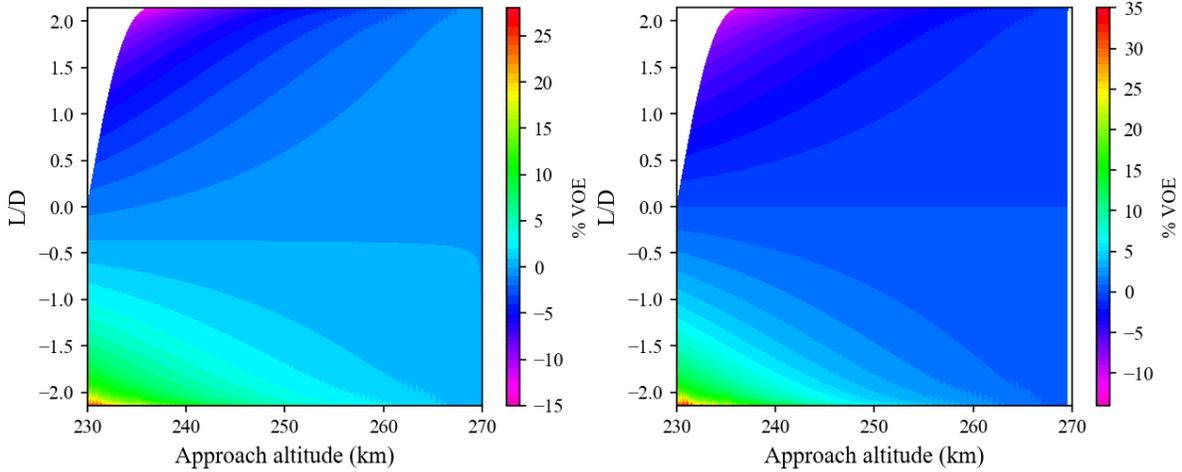

Figure 12. Contribution of VOE after approaching Venus AGAM (left) and PAGAM (right).

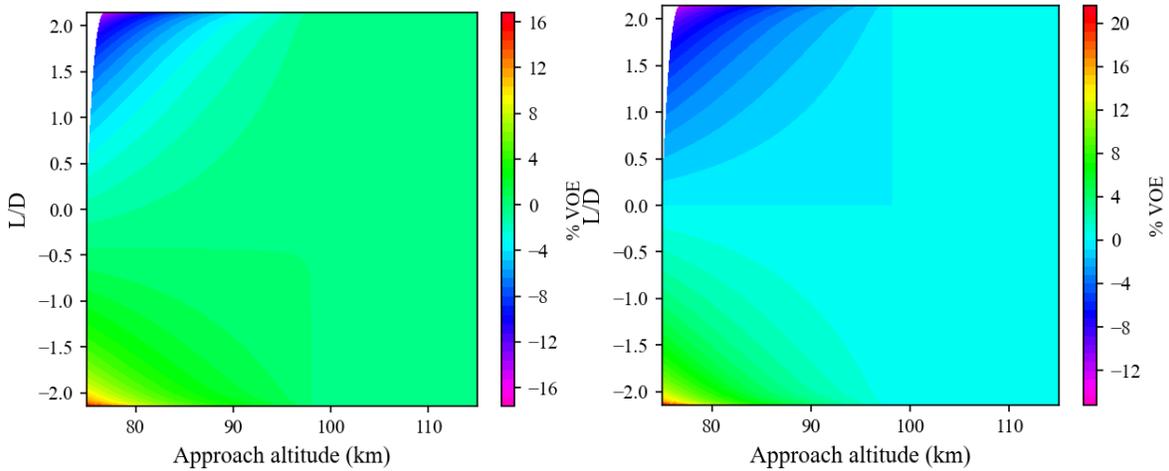

Figure 13. Contribution of VOE after approaching Mars AGAM (left) and PAGAM (right).

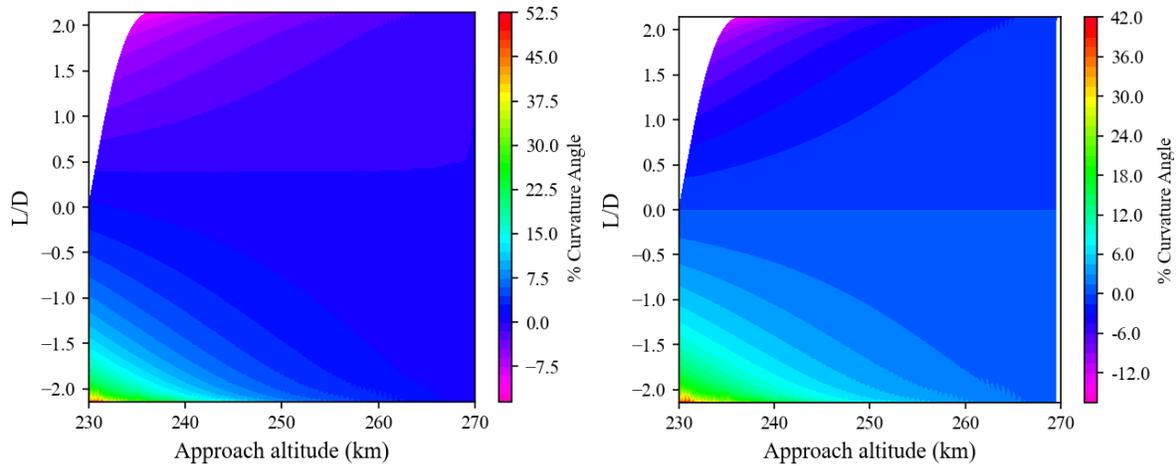

Figure 14. Contribution of curvature angle after approaching Venus AGAM (left) and PAGAM (right).

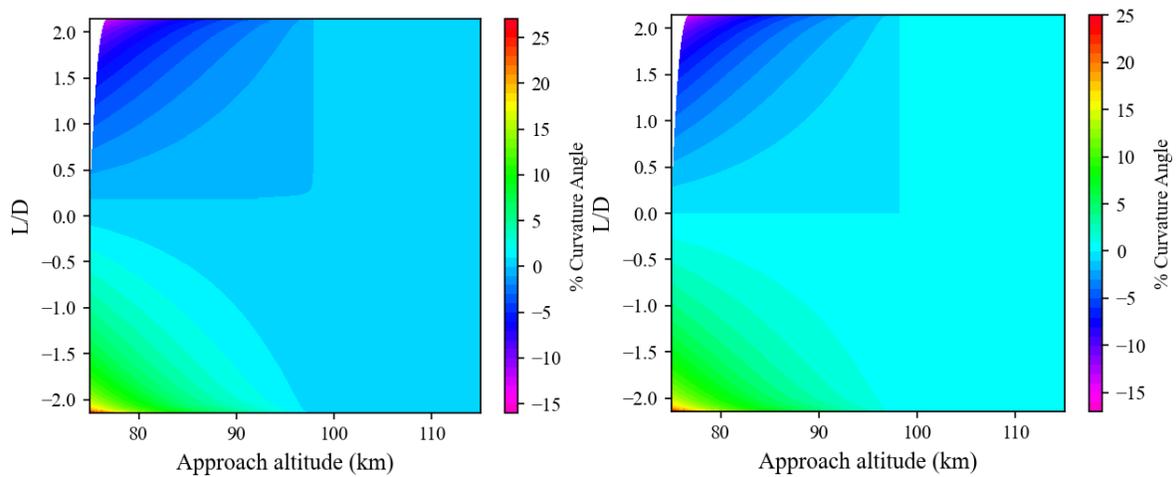

Figure 15. Contribution of curvature angle after approaching Mars AGAM (left) and PAGAM (right).

The PAGAM at Venus increases the energy gains by more than 7% for the closest approach with lower values of L/D. It is also observed a large region between L/D -0.5 to 0.5, where the effects of the aerodynamic forces present a contribution lower than ±5%, then, for this region, the GAM is closer to the AGAM and PAGAM because the Lift is lower than Drag. The same effect is observed in the largest altitudes, due to the lower density. In the case of Mars, the contributions of the maneuvers are equivalent, reducing the magnitude (see Figures 12 and 13). One application of the PAGAM without Lift is to increase the performance of the GAM, taking advantage of the proximity of the planet, compared to the GAM above the atmospheric boundary.

For the curvature angle, the AGAM and PAGAM at Venus, increase the angle by more than 40% when is compared to the GAM, in regions where L/D is minimum at high values of density (see Figure 14). The AGAM increases the curvature angle by more than 10% for Venus (figure 14) and 6% for Mars (figure 15), when compared to the PAGAM. Then, the increase in velocity magnitude of the PAGAM reduces the Time Of Flight and the influence of Drag required to increase the angle.

The results show a significant improvement with respect to the GAM, in terms of energy changes and turn angle, when the AGAM and/or the PAGAM are applied at upper atmospheric flight, which are feasible maneuvers. The effects reported here, are lower than the ones presented by Armellin et al. (2006); Qi and Ruiter (2021), which implemented optimal control to the maneuvers applied at Mars, but at the lowest altitudes (below 60 km), which include the limitations described in section 2.1.

Due to the use of continuum impulse for the PAGAM, is calculated the Delta-V required for the propulsion system to generate the Thrust with the same magnitude and in opposite direction to Drag, as a function of the altitude and the time of the atmospheric flight. The results are presented in Figure 16.

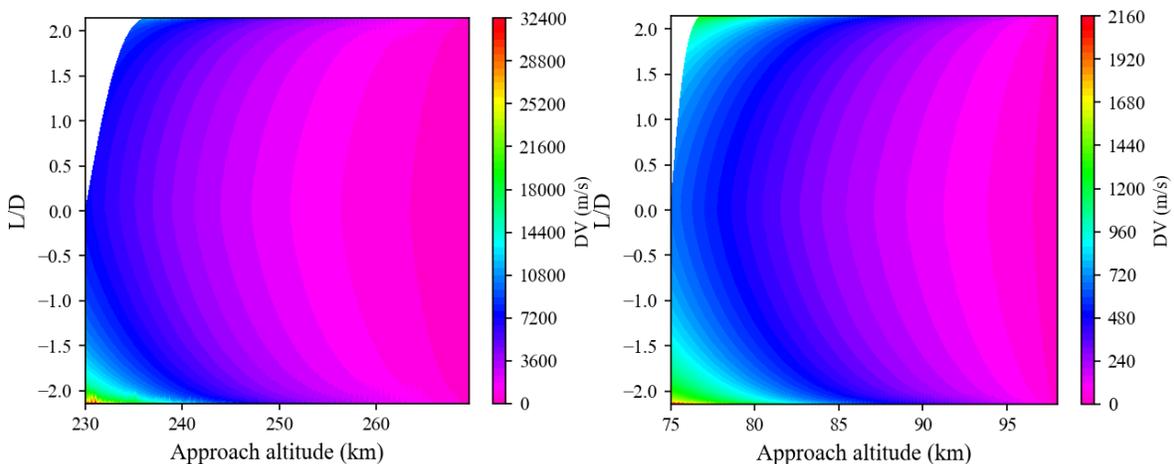

Figure 16. DV requires for the PAGAM Venus (left) and Mars (right).

The maximum Delta-V required is located at the same regions of the maximum VOE; minimum altitude and L/D ratios. The transition from blue to violet shows regions where the DV is lower, at larger altitudes, and without significant influence of the L/D due to the low values of density. However, the values for Venus show that the PAGAM requires the largest Delta-V. In other words, it is an expensive maneuver compared to the AGAM. In the case of Mars, the Delta-V values are achievable using conventional propulsion systems, showing a low cost compared to the PAGAM at Venus (Figure 16).

## 5. Conclusions

In this paper, it was presented a list of the closest approaches of historical missions applying the Gravity Assist Maneuver, showing the passage in the atmosphere. It was also presented the maneuvers derived from the Gravity- Assist, including impulse and under the influence of aerodynamic forces, with a detailed description of the Technological Readiness Levels and limitations for the application of the Aero-gravity Assisted Maneuvers and Powered Aero-gravity Assisted Maneuvers at lower altitudes or lower atmosphere.

For the maneuver in the upper atmosphere, it was possible to observe significant effects. The larger influence in the Variations of Energy was presented when the Lift was directed to the planet, at the lowest values of Lift-to-Drag ratio and close approach, increasing at the same time the variations in the angle of curvature due to the duration of the atmospheric flight. When the Powered Aero-gravity Assisted Maneuvers was applied, it was observed a significant increase in the Variation of Energy, when compared to the Aero-gravity Assisted Maneuver, due to the effect of reducing the losses by Drag, but at the same time presenting a small reduction in the angle of curvature. Due to the range of the altitude, the atmospheric properties, and the differences in the masses of the planet, it is observed larger differences between the maneuvers at Venus and Mars, showing a better performance in the case of Venus, where it is recommended the application of the Aero-gravity Assisted Maneuver; while for Mars, the Powered Aero-gravity Assisted Maneuver is the best choice due to the feasibility and low cost when is compared to Venus.

The Aero-gravity Assisted Maneuver and Powered Aero-gravity Assisted Maneuver present a deviation in the approach angle below 2% and altitude below 4% from the designed values. These results show that the maneuver presents significant changes in the energy gain and the angle of curvature compared to the Gravity Assisted Maneuver, without a significant deviation from the pericenter. At the same time, this allowed to verify the use of Broucke´s methodology from Gravity-Assist to Aero-gravity-Assisted Maneuver and Powered Aero-gravity Assisted Maneuver, showing good agreement with the desired values.

From the results and discussion, it is possible to conclude that the Powered and unpowered Aero-gravity Assisted Maneuvers at high altitudes increase the effects of the traditional Gravity Assisted Maneuver and they are technologically feasible. The use of the Aero-gravity Assisted Maneuver is suggested for increasing the gains in curvature and the Powered Aero-gravity Assisted Maneuver for increasing the gains in energy.

**Acknowledgments**
The authors wish to express their appreciation for the support provided by Science and Technology Institute (ICT – UNIFESP) and by the National Institute for Space Research (INPE).
This work was supported by grants # 2019/26605-2 and # 2016/24561-0, from São Paulo Research Foundation (FAPESP); grants # 406841/2016-0, 301338/2016-7 and 303102/2019-5 from the National Council for Scientific and Technological

Development (CNPq); grant # 88882.317514/2013-01 from the National Council for the Improvement of Higher Education (CAPES).

**List of figures**